\documentclass[3p]{elsarticle}
\usepackage{amsmath}
\usepackage{amsfonts,amssymb}

\def\ghol{\mathfrak {hol}}

\def\gsim{\mathfrak {sim}}
\def\gso{\mathfrak{so}}
\def\Hol{{\rm Hol}}
\def\Sim{{\rm Sim}}
\def\SO{{\rm SO}}
\newtheorem{Th}{Theorem}
\newtheorem{Prop}{Proposition}
\newtheorem{Cor}{Corollary}
\newtheorem{Lem}{Lemma}
\newtheorem{Def}{Definition}
\newtheorem{Ex}{Exercise }
\def\bt{\begin{Th}}
\def\et{\end{Th}}
\def\bp{\begin{Prop}}
\def\ep{\end{Prop}}
\def\bc{\begin{Cor}}
\def\ec{\end{Cor}}
\def\bl{\begin{Lem}}
\def\el{\end{Lem}}
\def\bd{\begin{Def}}
\def\ed{\end{Def}}
\def\bex{\begin{Ex}}
\def\eex{\end{Ex}}
\def\pf{\noindent{\it Proof:\ }}

\def\be{\begin{equation}}
\def\ee{\end{equation}}
\def\ben{\begin{enumerate}}
 \def\een{\end{enumerate}}
\def\ba{\begin{array}{rlll}}
\def\ea{\end{array}}
\def\bea{\begin{eqnarray}}
\def\eea{\end{eqnarray}}
\def\bean{\begin{eqnarray*}}
\def\eean{\end{eqnarray*}}

\catcode`@=11 \@addtoreset{equation}{section} \catcode`@=12

\def\Real{\mathbb{R}}
\def\g{\mathfrak{g}}
\def\h{\mathfrak{h}}

\def\Ric{\mathop{{\rm Ric}}\nolimits}
\def\p{\partial}
\def\R{\mathcal{R}}

\begin{document}

\begin{frontmatter}

\title{{Two-symmetric Lorentzian manifolds}}

\author{Dmitri V. Alekseevsky}
\address{The University of Edinburgh, School of Mathematics, The King's Buildings, Mayfield Road,
Edinburgh, {D.Aleksee@ed.ac.uk}}

\author{Anton S. Galaev}
\address{Department of Mathematics and Statistics, Faculty of Science, Masaryk University,
Kotl\'a\v rsk\'a~2, 61137 Brno, Czech~Republic,
{galaev@math.muni.cz}}


\begin{abstract} We classify two-symmetric Lorentzian
manifolds using  methods of the theory of holonomy groups. These
manifolds are exhausted by a  special type of pp-waves and,
  like the symmetric Cahen-Wallach spaces, they have commutative holonomy.
\end{abstract}

\end{frontmatter}

\section{Introduction}

Symmetric pseudo-Riemannian manifolds constitute  an important
class of spaces. A direct generalization of these manifolds is
provided by the so-called $k$-symmetric pseudo-Riemannian spaces
$(M,g)$ satisfying the condition $$\nabla^k R=0,\quad \nabla^{k-1}
R\neq 0,$$ where $k\geq 1$ and $R$ is the curvature tensor of
$(M,g)$. For Riemannian manifolds, the condition $\nabla^k R=0$
implies $\nabla R=0$ \cite{Tanno72}. On the other hand, there
exist pseudo-Riemannian $k$-symmetric spaces with $k\geq 2$, see
e.g. \cite{Kai,Sen08,Sen10}.

The fundamental paper by J.M. Senovilla \cite{Sen08}  is devoted
to a detailed investigation of two-symmetric Lorentzian spaces. It
contains many interesting results about such manifolds and their
physical applications. In particular, it is proven there that any
two-symmetric Lorentzian space admits a parallel null vector
field. A classification of four-dimensional two-symmetric
Lorentzian spaces is obtained in the paper  \cite{Sen10}, in which
it is shown that these spaces  are some special pp-waves. The
result is based on the Petrov classification of the Weyl tensors.

In the present paper we generalize the result of \cite{Sen10} to
any dimension. The main result can be stated as follows.

\begin{Th}\label{ThMain} Let $(M,g)$ be a locally indecomposable Lorentzian manifold of dimension
$n+2$. Then $(M,g)$ is two-symmetric if and only if locally there
exist coordinates $v,x^1,...,x^n,u$ such that
$$g=2dvdu+\sum_{i=1}^n(dx^i)^2+(H_{ij} u+F_{ij})x^ix^j(du)^2,$$
where   $H_{ij}$ is a nonzero diagonal real matrix with the
diagonal elements $\lambda_1\leq\cdots\leq\lambda_n$, and $F_{ij}$
is a  symmetric real matrix.

Any other metric of this form isometric to $g$ is given by the
same $H_{ij}$ and by $\tilde F_{ij}=cH_{ij}+F_{kl}a^k_ia^l_j$,
where  $c\in\Real$ and $a^j_i$ is an orthogonal matrix such that
$H_{kl}a^k_ia^l_j=H_{ij}$.
\end{Th}

By the Wu Theorem \cite{Wu}, any Lorentzian manifold $(M,g)$ is
either locally indecomposable, or it is locally a product of a
Riemannian manifold $(M_1,g_1)$, and of a locally indecomposable
Lorentzian manifold $(M_2,g_2)$. The manifold $(M,g)$ is
two-symmetric if and only if $(M_1,g_1)$ is locally symmetric and
$(M_2,g_2)$ is two-symmetric. Consequently, Theorem \ref{ThMain}
provides the complete local classification of two-symmetric
Lorentzian manifolds.

For the proof of Theorem \ref{ThMain}, we use the methods of the
theory of holonomy groups.  The assumption that a Lorentzian
manifold $(M,g)$ is two-symmetric implies that the holonomy
algebra $\ghol_m$ of $(M,g)$ at a point $m\in M$ annihilates the
tensor $\nabla R_m\neq 0$. This cannot happen if the holonomy
algebra is the whole Lorentzian algebra $\gso(1,n+1)$.  Hence the
holonomy algebra must preserve a null line and is contained in the
similitude algebra,
$\ghol_m\subset\gsim(n)=(\Real\oplus\gso(n))+\Real^n,$ the maximal
Lie algebra with this property \cite{ESI}. It is sufficient   to
consider the following two cases: $\ghol_m=\h+\Real^n$, where
$\h\subset\gso(n)$ is an irreducible subalgebra, and
$\ghol_m=\Real^n$.

We prove that the first case is impossible: for this we calculate
$\nabla R$ and $\nabla\Ric$, and show that the Weyl conformal
tensor $W$ is parallel ($\nabla W=0$). Then, using the results of
A.~Derdzinski and W.~Roter \cite{DR77,DR09} and of \cite{Galconf},
we get a contradiction.

The second case corresponds to pp-waves. The condition $\nabla^2
R=0$ and simple computations allow us to find the coordinate form
of the metric.

\section{Holonomy groups of Lorentzian manifolds}
 We recall some basic  facts  about holonomy groups of  Lorentzian manifolds that can be found in
 \cite{ESI,Gal2,Leistner}.
Let $(M,g)$ be a   Lorentzian $d$-dimensional  manifold  and
$\Hol^0 (M) = \Hol^0 (M)_m$ its connected holonomy group  at a
point $m \in M$. It is a subgroup of the (connected) Lorentz group
$ \SO(V)^0$ where $V = T_mM$ is the tangent space and it is
determined by its Lie algebra $\ghol(M) \subset  \gso(V)$ which is
called the holonomy algebra of $M$.

  The  manifold $(M,g)$ is locally indecomposable (i.e. locally is not a direct product
  of two pseudo-Riemannian manifolds) if and only if  the  holonomy group $\Hol^0 (M)$
  (or the holonomy algebra
   $\ghol(M)$) is weakly irreducible, i.e. it  does not preserve  any  proper nondegenerate  subspace of $V$.
 Any  weakly irreducible  holonomy group $\Hol (M)$  different from  the Lorentz group $\SO(V)^0$ is a  subgroup
 of  the horospheric  group  $\SO(V)_{[p]}$, the    subgroup of $\SO^0(V)$ which preserves a null line
 $[p] = \Real p$.
 This  group is identified with  the  group $\Sim_{n} = \Real^* \cdot \SO_{n} \cdot \Real^{n},\, n= d-2$ of similarity
 transformations  of   the
 Euclidean  space
 $E=\Real^n$ as follows (see \cite[Sect. 2.3]{ESI}).
  The  Lorentzian group $\SO(V)^0$ acts   transitively on the  celestial  sphere $S^{n} = P V^0$
  (the space of null lines)  which is
  the projectivization of the null cone $V^0 \subset V$ with the
  stabilizer $\SO(V)_{[p]}$. The stabilizer  has an open orbit $ S^{n} \setminus [p]$
  which is identified  via the stereographic  projection  with the Euclidean  space $E$.
  The group $\SO(V)_{[p]}$ acts in $E$ as the full connected Lie group of similarity transformations. Having in mind
  this isomorphism, we will call the   group  $\SO(V)_{[p]}$  the similarity group  and denote it by $\Sim_{n}$.

  Using the metric $<.,.> = g_m$, we will identify
 the Lorentz Lie algebra
$\gso(V) \simeq \gso(1,n+1)$  with the  space  $\Lambda^2V$ of
bivectors.
 Then  the Lie algebra $\gsim_{n}$ of the  similarity group  can be  written as
  $$  \gsim_{n} =\gso(V)_{[p]} = \Real p\wedge q + p\wedge E + \gso(E) $$
 where $p,q$ are  isotropic vectors with $<p,q> =1$  which span 2-dimensional Minkowski
   subspace $U$ and $E = U^\perp$ is its  orthogonal complement. The
   commutative ideal $p \wedge E$ generates   the  commutative  normal  subgroup $T_E\subset \Sim_{n}$
    which acts on $E$  by   translations.  This group is called  the vector group.
     The one-dimensional subalgebra
  $\Real p \wedge q = \gso(U)$  generates the   maximal diagonal  subgroup  ${\rm\bf A}$  of  $\Sim_{n}$
   which is  the
    Lorentz  group  $\SO(U)^0 $  and   the maximal  compact subalgebra
     $\gso(E)$ generates  the group  $\SO(E)$ of orthogonal  transformations of $E$.
 The  above decomposition of  the Lie algebra $\gsim_{n}$ defines  the Iwasawa decomposition
 $$\Sim_{n} = K\cdot {\rm\bf A}\cdot N = \SO(E)\cdot \SO(U)^0 \cdot T_E $$
 of the group $\Sim_{n}$.
 The list of connected weakly irreducible  holonomy groups $\Hol^0 (M)$ of
  Lorentzian manifolds  is known, see \cite{Leistner,ESI}.
  Assume  for simplicity  that $\Hol^0(M)$ is  an algebraic group.
  Then it contains  the vector group $T_E$   and     has  one  of the  following    forms:\\
 (type I)   $\Hol^0(M) = K \cdot \SO(U)^0 \cdot T_E$ \\
 (type II)  $\Hol^0(M) = K  \cdot T_E,$
 where $K \subset \SO(E)$ is  a connected holonomy group  of  a Riemannian $n$-dimensional manifold, i.e.   a product of  the Lie groups
  from the Berger list:
  $\SO_k,$ ${\rm U}_k,$ ${\rm SU}_k,$ ${\rm Sp}_1 \cdot {\rm Sp}_k$, ${\rm Sp}_k$, ${\rm G}_2$,  ${\rm Spin}_7$
  and  the  isotropy  groups of irreducible  symmetric Riemannian  manifolds.

   If  the  holonomy  group is not algebraic, it is  obtained  from one of  the
    holonomy groups of type I or II
   by  some twisting (holonomy groups of type III  and IV).
 Note that all these  holonomy groups act  transitively on the Euclidean space
 $E = PV^0 \setminus [p]$~\cite{Gal2}.

The Lorentzian holonomy algebras $\mathfrak{g}\subset\gsim(n)$ are
the following :

\begin{tabular}{llll} (type I) & $\Real p\wedge
q+\h+p\wedge E$, & (type II) & $\h+p\wedge E$,\\
(type III) & $\{\varphi(A)p\wedge q+A|A\in\h\}+p\wedge E$,& (type
IV) & $\{A+p\wedge \psi(A) |A\in\h\}+p\wedge E_1$,\end{tabular}

where $\h\subset\gso(E)$ is a Riemannian holonomy algebra;
$\varphi:\h\to\Real$ is a non-zero linear map that is zero on the
commutant~$[\h,\h]$; for the last algebra $E=E_1\oplus E_2$ is an
orthogonal decomposition, $\h$ annihilates $E_2$, i.e. $\h\subset
\gso(E_1)$, and $\psi:\h\to E_2$ is a surjective linear map that
is zero on the commutant $[\h,\h]$. The subalgebra
$\h\subset\gso(E)$, i.e. the $\gso(E)$-projection of $\g$ is
called {\it the orthogonal part of} $\g$.

A locally indecomposable simply connected Lorentzian manifold
admits a parallel null vector field if and only if its holonomy
group is of type II or IV.

  \section{The  holonomy group of a two-symmetric Lorentzian manifold}
   \bd A pseudo-Riemannian manifold $(M,g)$ with the curvature tensor $R$  is called a $k$-symmetric space  if
   $$ \nabla^k R = 0,\quad   \nabla^{k-1}R \neq 0.  $$
   \ed
   So,
    one-symmetric  spaces  are  the same as nonflat locally  symmetric   spaces ($\nabla R =0$, $R\neq 0$).
     Recall that  a  complete simply connected  locally symmetric space  is a  symmetric
     space,
     that is it admits   a  central symmetry $S_m$  with  center at any point $m$, i.e.  an involutive
      isometry $S_m$ which has $m$ as an isolated fixed  point.

Remark that for a Riemannian manifold the condition $\nabla^kR=0$
implies $\nabla R=0$ \cite{Tanno72}.

    All indecomposable  simply connected Lorentzian symmetric spaces are exhausted by
     the De Sitter  and the anti De Sitter spaces and
    the  Cahen-Wallach spaces, which have  the vector  holonomy  group $T_E$.

The following result is proven, using  so-called casual tensors
and the super-energy techniques,  in \cite{Sen08}.

      \bt\label{ThSenov}\cite{Sen08} Any two-symmetric Lorentzian manifold $(M,g)$  admits a parallel null
vector field.\et

This implies that the holonomy group can be only of type II or IV.
To make the exposition complete,  we will sketch a proof  of
Theorem \ref{ThSenov} using the holonomy theory.

The corner stone of the paper is the following statement.

  \bt\label{Thhol}  The  holonomy  group $\Hol^0(M)$ of an $(n+2)$-dimensional
  locally indecomposable two-symmetric Lorentzian manifold $(M,g)$  is the vector  group
  $T_E$ with the Lie algebra $p \wedge E \subset \gso(V)$. \et

It is known that any $(n+2)$-dimensional  Lorentzian manifold with
the holonomy algebra $p\wedge E$ is a pp-wave (see e.g. \cite[Sect
5.4]{ESI}), i.e. locally there exist coordinates $v,x^1,...,x^n,u$
such that the metric $g$ can be written in the form
$$g=2dvdu+\delta_{ij}dx^idx^j+H(du)^2,\quad \p_v H=0.$$
We will need only to decide which functions $H$ correspond to
two-symmetric spaces.

\subsection{Algebraic curvature tensors}

For a subalgebra $\g\subset\gso(V)$ define {\it the space of
algebraic curvature tensors of type} $\mathfrak{g}$,
$$\R(\mathfrak{g})=\{R\in \Lambda^2 V^*\otimes\mathfrak{g}\,|\, R(u,v)w+R(v, w)u+R(w, u)v=0
\text{ for all } u,v,w\in V\}.$$  If $\mathfrak{g}\subset\gso(V)$
is the holonomy algebra of   a  manifold $(M,g)$, where $V=T_mM$
is tangent space at some point $m\in M$,
 then the curvature tensor $R_m$ of $(M,g)$ belongs to
$\R(\mathfrak{g})$. The spaces $\R(\g)$ for holonomy algebras of
Lorentzian manifolds are found in \cite{Gal1,onecomp}. For
example, let $\g=\Real p\wedge q+\h+p\wedge E$. For a subalgebra
$\h\subset\gso(n)$ define the space
$$\mathcal{P}(\h)=\{P\in E^*\otimes\h\,|\,g(P(x)y,z)+g(P(y)z,x)+g(P(z)x,y)=0\text{ for all }
x,y,z\in E\}.$$ Any $R\in\R(\g)$ is uniquely determined  by the
data $(\lambda,e,P,R^0,T)$, where $$ \lambda\in\mathbb{R},\ e\in
E,\ P\in\mathcal{P}(\mathfrak{h}),\
R^0\in\mathcal{R}(\mathfrak{h}),\ T\in S^2E,$$ i.e. $T$ is a
symmetric tensor considered as an endomorphism of E. The tensor
$R$ is defined by
\begin{align*}
R(p,q)=& -\lambda p\wedge q-p\wedge e,\qquad R(X,Y)=R^0(X,Y)-p\wedge (P(Y)X-P(X)Y),\\
R(X,q)=&-g(e,X)p\wedge q+P(X)-p\wedge T(X),\qquad R(p,X)=0,\quad
\forall X,Y\in E.
\end{align*}   We will write
$$R=R^{(\lambda,e,P,R^0,T)}.$$ If some of these elements are zero,
we omit them. For example, if $R$ is defined only by $T$, then we
write $R=R^T$. Note that $$R^T=\sum_{i,j}T_{ij}p\wedge e_i\vee
p\wedge e_j,\quad T_{ij}=g(Te_i,e_j),\qquad\R(p\wedge E)=\{R^T|\,
T\in S^2 E\}\simeq S^2 E,$$ where $e_1,...,e_n$ is an orthonormal
basis of $E$, and $\vee$ denotes the symmetric product. Similarly,
$$\R(\h+p\wedge E)=\{R^{(P,R^0,T)}|\, P\in{\mathcal P}(\h),\,
R^0\in\R(\h),\,T\in S^2E\}.$$

 Now we define the space of covariant derivatives of the curvature
 tensor
$$\nabla\R(\mathfrak{g})=\{S\in {\rm Hom}(V,\R(\mathfrak{g}))=V^*\otimes\R(\mathfrak{g})\,|\, S_u(v,w)+S_v(w,u)+S_w(u,v)=0
\text{ for all } u,v,w\in V\}.$$ If $\mathfrak{g}\subset\gso(V)$
is the holonomy algebra of a  manifold $(M,g)$ at a point $m\in
M$, then $\nabla R_m\in\nabla\R(\mathfrak{g})$. The decomposition
of the space $\nabla\R(\gso(r,s))$ into irreducible
$\gso(r,s)$-modules is found in \cite{Str88}, see also
\cite{nablaR}.

It is not difficult to find the space $\nabla\R(\g)$ for each
Lorentzian holonomy algebra $\g\subset\gsim(n)$. It consists of
tensors $$S\in{\rm Hom}(V,\R(\g)),\quad S:u\in V\mapsto
S_u=R^{(\lambda_u,e_u,P_u, R^{0}_{u},T_u)}\in \R(\g)$$ satisfying
the second Bianchi identity. For example,
$$\nabla\R(p\wedge E)= \{S=q'\otimes R^{T}|\, T\in
S^2E\}\oplus\{S=R^{Q.}|\, Q\in S^3E\}\simeq S^2E\oplus S^3E,$$
here $q'=g(p,\cdot)$ is the 1-form $g$-dual to $p$, the tensor
$S=R^{Q.}$ is defined by $S_p=S_q=0$, $S_x=R^{Q_x}$, $x\in E$,
$Q_x\in S^2E$ (since $Q\in S^3E$).

\subsection{Adapted coordinates and  reduction lemma}

Let $(M,g)$ be an $(n+2)$-dimensional
  locally indecomposable (hence with weakly irreducible holonomy algebra $\g$)
  two-symmetric Lorentz manifold, i.e. the tensor $\nabla R$ is
  nonzero, parallel and annihilated by the holonomy algebra.
The space $\nabla\R(\gso(1,n+1))$ does not contain nonzero
elements annihilated by $\gso(1,n+1)$, see e.g. \cite{Str88}.
Since $\gso(1,n+1)$ is the only  irreducible holonomy algebra
\cite{ESI}, it follows that $\g\subset\gsim(n)$.

Let $(M,g)$ be a Lorentzian manifold  with the holonomy algebra
$\g\subset\gsim(n)$. Then $(M,g)$ admits a parallel distribution
of null lines.  According to \cite{Walker}, locally there exist so
called Walker coordinates $v,x^1,...,x^n,u$ such that the metric
$g$ has the form
\begin{equation}\label{Walker} g=2dvdu+h+2Adu+H (d
u)^2,\end{equation} where $h=h_{ij}(x^1,...,x^n,u)d x^id x^j$ is
an $u$-dependent family of Riemannian metrics,  $A=A_i(x^1,
\ldots, x^n,u)d x^i$ is an $u$-dependent family of one-forms, and
$H=H(v,x^1,...,x^n,u)$ is a local function on $M$. Consider the
local frame
$$p=\p_v,\quad X_i=\p_i-A_i\p_v,\quad q=\p_u-\frac{1}{2}H\p_v.$$
Let $E$ be the distribution generated by the vector fields
$X_1$,...,$X_n$. Clearly, the vector fields $p$, $q$ are
isotropic, $g(p,q)=1$, the restriction of $g$ to $E$ is positive
definite, and $E$ is orthogonal to $p$ and $q$. The vector field
$p$ defines the parallel distribution of null lines and it is
recurrent, i.e. $\nabla p=\theta\otimes p$, where
$\theta=\frac{1}{2}\p_vHdu$. Since the manifold is locally
indecomposable, any other recurrent vector field is proportional
to $p$. Next, $p$ is proportional to a parallel vector field if
and only if $d\theta =0$, which is equivalent to
$\p^2_vH=\p_i\p_vH=0$. In the last case the coordinates can be
chosen in such a way that $\p_v H=0$ and $\nabla p=\nabla \p_v=0$,
see e.g. \cite{ESI}.

Let $\g\subset\gsim(n)$ be the holonomy algebra of the Lorentzian
manifold $(M,g)$ and $\h\subset\gso(E)$ be its orthogonal part.
Then there exist the decompositions
\begin{equation}\label{LM0A}E=E_0\oplus E_1\oplus\cdots\oplus E_r,\quad
\h=\{0\}\oplus\h_1\oplus\cdots\oplus\h_r\end{equation} such that
$\h$ annihilates $E_0$, $\h_i(E_j)=0$ for $i\neq j$, and
$\h_i\subset\gso(E_i)$ is an irreducible subalgebra for $1\leq
i\leq s$.  Ch.~Boubel \cite{Boubel} proved  that there exist
Walker coordinates
$$v,\,x_0=(x_0^1,\ldots ,x_0^{n_0}),\ldots ,x_r=(x_r^1,...,x_r^{n_r}),\,u$$
adapted to the decomposition \eqref{LM0A}. This means that \be
\label{h=sum}h=h_0+h_1+\cdots+h_r,\quad h_0=\sum_{i=1}^{n_0}(d
x_0^i)^2,\quad h_\alpha=\sum_{i,j=1}^{n_\alpha}h_{\alpha ij} d
x_{\alpha}^id x_{\alpha}^j,\ee
$$A=\sum_{\alpha=1}^r A_\alpha,\quad A_0=0,\quad A_\alpha=\sum_{k=1}^{n_\alpha} A^\alpha_k d x^k_\alpha,$$
and   one has \be\label{halpha=} \frac{\p}{\p
{x^k_\beta}}h_{\alpha ij}=\frac{\p}{\p
{x^k_\beta}}A^\alpha_i=0,\quad  \text{ if }\beta\neq\alpha.\ee
 The coordinates can be chosen
 so that $A=0$, see \cite{GL10}. Thus we will assume that $g$ is given by
 \eqref{Walker} with $A=0$, and with $h$ satisfying \eqref{h=sum}
 and \eqref{halpha=}.

For $\alpha=0,...,r$, consider  the submanifolds
$M_{\alpha}\subset M$ defined by $x_\beta=c_\beta$,
$\alpha\neq\beta$, where $c_\beta$ are constant vectors. Then the
induced  metric is given by
$$g_\alpha=2dvdu+h_\alpha+H_\alpha(du)^2.$$

\begin{Lem}\label{Lemgi} The submanifold $M_\alpha\subset M$ is
totally geodesic. The orthogonal part
 of the holonomy algebra
$\g_{\alpha}$ of the metrics $g_\alpha$ coincides with
$\h_\alpha\subset\gso(E_\alpha)$, which is irreducible for
$\alpha=1,...,r$. If the metric $g$ is two-symmetric, then the
curvature tensor of each metric $g_\alpha$ satisfies
$\nabla^2R=0$.
 \end{Lem}

{\it Proof.} The non-zero Christoffel symbols of the metric
\eqref{Walker} with $A=0$ are the following:
$$\Gamma^v_{uu}=\frac{1}{2}H_{,u},\quad
\Gamma^v_{iu}=\frac{1}{2}H_{,i},\quad
\Gamma^v_{vu}=\frac{1}{2}H_{,v},\quad
\Gamma^i_{ju}=\frac{1}{2}h^{ik}h_{jk,u},$$
$$\Gamma^i_{uu}=-\frac{1}{2}h^{ik}H_{,k},\quad
\Gamma^u_{uu}=\frac{1}{2}H_{,v},\quad
\Gamma^i_{jk}=\Gamma^i_{jk}(h),$$ where the comma denotes the
partial derivative and $\Gamma^i_{jk}(h)$ are the Christoffel
symbols of the metric $h$. This shows that the Christoffel symbols
of the metric $g_\alpha$ are equal to the restrictions to
$M_\alpha$ of the corresponding Christoffel symbols of the metric
$g$, i.e. each submanifold $M_\alpha\subset M$ is totally
geodesic. This implies that if $\nabla^2 R=0$, then each
$g_\alpha$ satisfies the same condition. Finally, the statement
about the orthogonal parts follows from the fact that the
orthogonal part of any Walker metric $g$ coincides with the
holonomy algebra of the induced connection on the vector bundle
with the fibers $p_m^\bot/\Real p_m\simeq E_m$, and this
connection does not depend on the function $H$ \cite{ESI}. \qed

\subsection{Sketch of the proof of Theorem \ref{ThSenov} using the holonomy theory}

 We may assume that the metric $g$ is locally given by \eqref{Walker} with $A=0$,
 and with $h$ satisfying \eqref{h=sum}
 and \eqref{halpha=}. As it is noted above, it is enough to prove that
$\p^2_vH=\p_i\p_vH=0$. Clearly, this will be true if it is true
for each metric $g_\alpha$.

\bl If $\g$ is the holonomy algebra of type I (with any orthogonal
part $\h\subset\gso(E)$), or $\g$ is the holonomy algebra of type
III with an irreducible orthogonal part $\h\subset\gso(E)$, then
the subspace $\nabla\R(\g)^0\subset\nabla\R(\g)$, consisting of
tensors annihilated by $\g$, is trivial. \el

{\bf Proof.} If $\g$ is of type I, then it contains $A=p\wedge q$.
If $\g$ is of type III, then $\h\subset
\mathfrak{u}(E)\subset\gso(E)$ and for some $a\in\Real$, the
element $A=p\wedge q+aJ$ belongs to $\g$. The lemma follows from
the consideration of the tensors in $\nabla\R(\g)$ annihilated by
 the operator $A$ and the second Bianchi identity  as in Lemma
\ref{Lem2} below. \qed

The lemma shows that the holonomy algebra of each metric
$g_\alpha$ cannot be of type I or III, i.e. it is of type II or
IV. Thus, $\p^2_vH=\p_i\p_vH=0$ holds. \qed

\subsection{Proof of Theorem \ref{Thhol}}

Consider the decomposition \eqref{LM0A}. If $E=E_0$, then $\h=0$
and there is nothing to prove. If $E_1\neq 0$, then the metric
$g_1$ satisfies $\nabla^2 R=0$ and the orthogonal part of its
holonomy algebra  $\h_1\subset\gso(E_1)$ is irreducible. We will
show that this is not possible.

Thus, suppose that $(M,g)$ satisfies $\nabla^2R=0$ and its
holonomy algebra equals to  $\g=\h+p\wedge E$, where
$\h\subset\gso(E)$ is irreducible.

\bl\label{Lem2} Let $\g=\h+p\wedge E$, where $\h\subset\gso(E)$ is
irreducible. Then the subspace $\nabla\R(\g)^0\subset\nabla\R(\g)$
of $\g$-annihilated tensors is the one-dimensional subspace given
by
$$\nabla\R(\g)^0=\Real S,\quad S=q'\otimes R^{{\rm
Id}_E},\quad q'=g(p,\cdot).$$ \el

{\it Proof.} Let $S\in\nabla\R(\g)^0$. For any $v\in V$, the
element $S_v\in\R(\g)$ can be written as $S_v=R^{(R^{0}_{v},
P_v,T_v)}$ as it is explained above. Since $S_\cdot(p,\cdot)=0$,
by the second Bianchi identity $S_p=0$. The fact that $\g$
annihilates $S$ can be expressed as
$$[A,S_{v_3}(v_1,v_2)]-S_{A v_3}(v_1,v_2)-S_{v_3}(A
v_1,v_2)-S_{v_3}(v_1,A v_2)=0$$ for all $A\in\g$ and
$v_1,v_2,v_3\in V$. Let $U,X,Y,Z\in E$. We have
$$[p\wedge X,S_U(Y,Z)]=0.$$
Hence, $R^{0}_{U}(Y,Z)X=0$, i.e. $R^{0}_{U}=0$. Next, $$[p\wedge
X,S_Z(Y,q)]-S_Z(Y,X)=0.$$ Consequently,
$$-p\wedge P_Z(Y)X-p\wedge (P_Z(Y)X-P_Z(X)Y)=0,$$
i.e. $2P_Z(Y)X=P_Z(X)Y$. Since this equality holds for any $X,Y\in
E$, we conclude $P_Z=0$. We have got $S_Z(X,Y)=0$. Similarly,
$$[p\wedge X,S_q(Y,Z)]=0,$$ i.e. $R^{0}_{q}=0$.
The equality $$[p\wedge X,S_q(Y,q)]-S_X(Y,q)-S_q(Y,X)=0$$ implies
$$T_X(Y)=2P_q(Y)X-P_q(X)Y.$$
From the second Bianchi identity $$S_q(X,Y)+S_X(Y,q)+S_Y(q,X)=0$$
it follows that $$T_X(Y)-T_Y(X)=P_q(X)Y-P_q(Y)X.$$ We conclude
$P_q(Y)X-P_q(X)Y=0$. This and the definition of the space
${\mathcal P}(\h)$ imply $P_q=0$. Consequently, $T_X=0$. Finally,
let $A\in\h$, then
$$[A,S_q(X,q)]-S_q(AX,q)=0.$$ This implies $AT_q(X)=T_q(AX),$ i.e. $T_q$
commutes with $\h$. Since $T_q$ is a symmetric endomorphism of $E$
and $\h\subset\gso(E)$ is irreducible, by the Schur Lemma, $T_q$
is proportional to the identity. This proves the lemma. $\Box$

We  write the metric $g$ in the form \eqref{Walker}. Then  $\p_v$
is parallel and $\p_vH=0$.

By Lemma \ref{Lem2},  $\nabla R$ has the form
\begin{equation}\label{nabR}\nabla_U R=fg(p,U)R^{{\rm Id}_E},\quad \forall U\in TM,\end{equation}
for some smooth function $f$. It is clear that
$$R^{{\rm Id}_E}(U_1,U_2)=p\wedge ((U_1\wedge U_2)p),\quad \forall U_1,U_2\in TM.$$

\bl Under the above assumptions, the conformal Weyl curvature
tensor $W$ is parallel, i.e.~$\nabla W=0$. \el

{\bf Proof.} It is known that $$W=R+L\wedge g,$$ where
$$L=\frac{1}{d-2}\left(\Ric-\frac{s}{2(d-1)}{\rm
Id}\right)$$ is the Schouten tensor, $\Ric$ is the Ricci operator,
and $s$ is the scalar curvature. Recall that by definition,
$$(L\wedge g)(U_1,U_2)=LU_1\wedge U_2+U_1\wedge LU_2,\quad U_1,U_2\in TM.$$
 For any vector field $U$ it
holds $$\nabla_UW=\nabla_UR+(\nabla_UL)\wedge g.$$ Let the indexes
$a,b$ run from $0$ to $n+1$, and let $X_0=p$, $X_{n+1}=q$. The
covariant derivative of the Ricci operator is given by
\begin{multline*}(\nabla_{U_1}{\rm
Ric})U_2=g^{ab}\nabla_{U_1}R(U_2,X_a)X_b=g^{ab}fg(p,U_1)R^{{\rm
Id}_E}(U_2,X_a)X_b= g^{ab}fg(p,U_1)(p\wedge((U_2\wedge
X_a)p))X_b\\=fg(p,U_1)(g^{ab}g(p,X_b)(U_2\wedge
X_a)p-g^{ab}g((U_2\wedge X_a)p,X_b)p)\\=fg(p,U_1)((U_2\wedge
p)p-g^{ab}g(g(U_2,p)X_a-g(X_a,p)U_2,X_b)p)=(2-d)fg(p,U_1)g(p,U_2)p.\end{multline*}
Thus, $\label{nablaRic} (\nabla_{U_1}{\rm
Ric})U_2=-nfg(p,U_1)g(p,U_2)p.$ The gradient of the scalar
curvature is given by $$g({\rm
grad}s,U_1)=g^{ab}g((\nabla_{U_1}{\rm Ric})X_a,X_b)=0,$$ i.e.
${\rm grad}s=0$. Hence,
$$(\nabla_{U_1}L)U_2=-fg(p,U_1)g(p,U_2)p.$$
Consequently, \begin{multline*}(\nabla_{U_1}L)U_2\wedge
U_3+U_2\wedge (\nabla_{U_1}L)U_3=-fg(p,U_1)g(p,U_2)p\wedge
U_3-U_2\wedge fg(p,U_1)g(p,U_3)p\\=fg(p,U_1)(p\wedge
g(p,U_3)U_2-p\wedge g(p,U_2)U_3)=-fg(p,U_1)p\wedge ((U_2\wedge
U_3)p)=-\nabla_{U_1}R(U_2,U_3).\end{multline*} Thus,
$(\nabla_{U_1}L)\wedge g=-\nabla_{U_1}R$ and  $\nabla W=0.$ \qed

The condition $\nabla W=0$ under the above assumptions implies
that $(M,g)$ is a pp-wave. If $W=0$, then this is proved in
\cite{Galconf}. If $W\neq 0$, then the results of  A.~Derdzinski
and W.~Roter  \cite{DR09,DR77} show that either $\nabla R=0$, or
$(M,g)$ is a pp-wave.
 Thus the holonomy algebra of $(M,g)$ is contained in $p\wedge E$ and we get
a contradiction. This proves the theorem. \qed

\section{Lorentzian manifolds with vector  holonomy group $T_E$ (pp-waves)}

   In this section we  derive   formulas for  the  curvature tensor and
    its covariant derivatives    for  an $(n+2)$-dimensional  Lorentzian
     manifold  with  the vector holonomy group $\Hol(M) = T_E$
   (or, equivalently, the  holonomy algebra $\ghol(M) = p \wedge E$).

 \subsection{  Adapted local  coordinates and associated  pseudo-group of  transformations }

 It is well known  that
the connected holonomy group of a Lorentzian manifold $(M,g)$ is a
subgroup of $T_E$ if and only if in  a neighborhood of any point
$x \in M$
  with  respect to some local coordinates $v, x^1, \cdots , x^n, u$
   (called adapted coordinates) the metric is given by
 \be\label{pp-wave}   g = 2 du dv + \delta_{ij}dx^i dx^j  + Hdu^2,
 \ee
where $H$ is a function  of $x^i$ and $u$, see e.g. \cite[Sect.
5.4]{ESI}. Such Lorentzian manifolds are called pp-waves.

\bl Any two adapted coordinate systems with the same $\partial_v$
 are
related by \be \label{transform} \tilde v = v  -\sum_j
a^j_i\frac{d b^j(u)}{d u}x^i+d(u), \quad  \tilde x^i =
a^i_jx^j+b^i(u), \quad \tilde u = u+c, \ee
 where  $c\in\Real$, $a^j_i$ is an orthogonal matrix, and $b^i(u)$, $d(u)$ are arbitrary functions of
 $u$.\el

 {\bf Proof.} In \cite{GL10} it is shown that two Walker
systems of coordinates with the same $\p_v$ are related by
$$\tilde
 v=v+f(x^1,\dots,x^n,u),\quad \tilde
 x^i=\phi^i(x^1,\dots,x^n,u),\quad \tilde u=u+c.$$
Since $h=\delta_{ij}dx^i dx^j$ must be preserved,
$\phi^i(x^1,\dots,x^n,u)$ must define an $u$-dependent family of
isometries of $\Real^n$, i.e. $$\phi^i(x^1,\dots,x^n,u)=
a^i_j(u)x^j+b^i(u),$$ where $a^j_i(u)$ is a family of orthogonal
matrices. Next, the equalities
$g(\p_i,\p_u)=g(\tilde\p_i,\tilde\p_u)=0$ imply $$\p_i f+\sum_k
a^k_i(u)\frac{d}{du}\left(a^k_r(u)x^r+b^k(u)\right)=0.$$ This
shows that $\sum_k a^k_i(u)\frac{d}{du}a^k_r(u)=0$, i.e.
$\frac{d}{du}a^k_r(u)=0$. Finally, we  easily find the function
$f$. \qed

 \subsection{Levi-Civita  connection}
 We associate  with  an  adapted  coordinate system  $(u, x^i, v)$ of  a pp-wave  space $(M,g)$
 with a potential  $H= H(x^i,u)$
 a   standard   field of frames
 $$ p= \p_v,\quad  e_i = \p_i,\quad q = \p_u - \frac{1}{2}H\p_v $$
 and the dual field of coframes
 $$ p' = dv + \frac{1}{2}H du, \quad e^i = dx^i,\quad q' = du. $$
  The Gram matrix  of these bases is given by

 $$
 G  = \left( \begin{array}{ccc}
  0 & 0 & 1\\
 0 & \b1_n &0 \\
 1 & 0& 0
 \end{array}
 \right).
 $$

  We will  consider coordinates of all  tensor fields with respect to these non-holonomic frame and coframe.
   Then the covariant  derivative
  of a vector $Y = Y^p p + Y^i e_i + Y^q  q$ and a covector
  $\omega = \omega_p p' + \omega_i e^i + \omega_q q'$
  in direction of a vector field $X$ can be written  as
  $$\nabla_XY = \p_X Y + A_XY , \,\,  \nabla_X \omega = \p_X \omega - A_X^T \omega $$
 where $\p_X$ is the coordinate derivative     in direction of  $X$  and $A_X$ is  a matrix and
 $A_X^T $ is the transposed  matrix.
 \bl The  matrices $A_u, A_i, A_v$ of  the  connection   which correspond to the  coordinate vector fields
 $\p_u, \p_i, \p_v$  and their transposes   are given by

$$ A_u
  = \left( \begin{array}{ccc}
  0 & \frac{1}{2}H_{,i} & 0\\
 0 & 0 &-\frac{1}{2}H_{,i} \\
 0 & 0& 0
 \end{array}
 \right),
\,\,\,
A^T_u
  = \left( \begin{array}{ccc}
  0 & 0 & 0\\
 \frac{1}{2}H_{,i} & 0 &0 \\
 0 & -\frac{1}{2}H_{,i}& 0
 \end{array}
 \right),
\,\,\,
A_i = A_i^T = A_v = A_v^T = 0.
  $$
    In particular, $\nabla p  = \nabla p'=0$.

 \el
 \pf  The  only non zero Christoffel symbols  are
 $$\Gamma^v_{uu} = \frac{1}{2}H_{,u},\quad  \Gamma^i_{uu} = -\frac{1}{2} H_{,i} , \quad \Gamma^v_{iu}=
 \frac{1}{2}H_{,i} $$
  where the commas stand for   the partial derivatives.
  Then we calculate
  $$\nabla \p_v= \nabla p = 0,\quad \nabla_u \p_i = \frac{1}{2}H_{,i}p,
  \quad
   \nabla_u q= \nabla_u(\p_u - \frac{1}{2}H \p_v)= \frac{1}{2}H_{,u}p -\frac{1}{2} H_{,i} e_i -\frac{1}{2}
  H_{,u}p
   = - \frac{1}{2}H_{,i}e_i, $$
 $$ \nabla_i \p_j =0,\quad \nabla_i \p_u = \frac{1}{2}H_{,i} p,\quad
   \nabla_i q = \nabla_i (\p_u - \frac{1}{2}H p) =0,
  \quad
 \nabla_v \p_u = \nabla_v \p_i = \nabla_v \p_v =0.
 $$\qed

 \bc

 A Lorentzian  manifold $M$ with vector  holonomy group $\Hol(M) = T_E$ has  the
 (globally defined) parallel vector field $p = \p_v$  and  parallel 1-form  $q'=du$.

\ec
\subsection{The curvature  tensor of a pp-wave space}

 \bl With respect to the standard  frame $p= \p_v,\, e_i= \p_i,\, q=\p_u - \frac{1}{2}H \p_v$
  and  the dual coframe $p',e^i, q'$, the  curvature  tensor  of a pp-wave  with potential $H(u,x^i)$
 is given by
 $$  R=  \sum_{i,j} \frac{1}{2}H_{,ij}(p \wedge e_i \vee p \wedge e_j) \,\, {\text {( the contravariant curvature  tensor)} }
$$
 $$  \bar R=  \frac{1}{2}H_{,ij}(q' \wedge e^i \vee q' \wedge e^j ) \,\,\,  {\text {( the covariant curvature  tensor)}}.
$$
\el

\pf It follows from the  formula $R(X,Y) = \p_X A_Y - \p_Y A_X -
A_{[X,Y]}$. \qed \bc The Ricci tensor of  $M$ is given by
$${\rm ric} = \frac{1}{2}\Delta H q' \otimes q' = \frac{1}{2} \Delta H du^2
$$
where $\Delta$ is the Laplacian in $\Real^n$. \ec

\subsection{The covariant derivatives of the curvature tensor}
Note that for any $i,j$, the covariant tensor  $q' \wedge e^i \vee
q' \wedge e^j$   and the  contravariant tensor $p \wedge e_i \vee
p \wedge  e_j$ are parallel. Hence  the first covariant derivative
of the curvature tensor is given by \be \nabla \bar
R=\frac{1}{2}H_{,ijk}e^k\otimes(q' \wedge e^i \vee q' \wedge e^j
)+\frac{1}{2}H_{,iju}q'\otimes(q' \wedge e^i \vee q' \wedge e^j
).\ee

 \bc\label{Corsym} The manifold $(M,g)$ is a locally symmetric space
 if and only if  the Hessian $H_{,ij}$ of the potential $H$  is  a constant,
 that is $H = H_{ij}x^i x^j+G_i(u) x^i + K(u)$.
 \ec
It can be shown that in the last case the coordinates can be
chosen in such a way that
$H=\lambda_1(x^1)^2+\cdots+\lambda_n(x^n)^2$ for some non-zero
real numbers $\lambda_i$ such that
$\lambda_1\leq\cdots\leq\lambda_n$ \cite{C-W}.

The second covariant derivative of the curvature tensor is given
by
\begin{multline} \nabla^2 \bar
R=\left(\frac{1}{2}H_{,ijk}-\frac{1}{4}\sum_kH_{,k}H_{,ijk}\right)q'^2
\otimes (q' \wedge e^i \vee q' \wedge e^j)\\ +
 \frac{1}{2}H_{,ijku}( q' \vee e^k) \otimes (q' \wedge e^i \vee q' \wedge e^j) +
 \frac{1}{2}H_{,ijk\ell}(  e^k \otimes e^{\ell}) \otimes (q' \wedge e^i \vee q' \wedge
 e^j).\end{multline}

This implies the following.

 \bt\label{Th2}
 A pp-wave with the metric \eqref{pp-wave} is two-symmetric
if and only if $$H = (u H_{ij}+F_{ij})x^i x^j+ G_i(u) x^i +K(u),$$
where $H_{ij}$ and $F_{ij}$ are symmetric real matrices, the
matrix $H_{ij}$ is non-zero, $G_i(u)$ and $K(u)$ are functions of
$u$. \et

\section{Proof of Theorem \ref{ThMain}}
To prove the theorem we  start with the metric \eqref{pp-wave} and
$H$ as in Theorem \ref{Th2} and use transformation
\eqref{transform} in order to write the metric as in Theorem
\ref{ThMain}. Let $\tilde v,\tilde x^1,...,\tilde x^n,\tilde u$ be
a new coordinate system. We may assume that the inverse
transformation is given by \be\label{trans1} u =\tilde u+c, \quad
x^i = a^i_j\tilde x^j+b^i(\tilde u), \quad v = \tilde v  -\sum_j
a^j_i\frac{d b^j(\tilde u)}{d \tilde u}\tilde x^i+d(\tilde u).\ee
For the new function $\tilde H$ written as in Theorem \ref{Th2} we
get
\begin{align}\label{tHkl} \tilde
H_{kl}&= H_{ij}a^i_ka^j_l,\\\label{tFkl}
 \tilde
F_{kl}&=(cH_{ij}+F_{ij})a^i_ka^j_l,\\
\label{tGk}\tilde G_k(\tilde u)&=-2\sum_j a^j_k\frac{d^2 b^j}{(d
\tilde u)^2}+2((\tilde u+c)H_{ij}+F_{ij})b^ia^j_k+G_ia^i_k,\\
\tilde K(\tilde u)&=2\frac{dd(\tilde u)}{d\tilde u}
+\sum_j\left(\frac{db^j}{du}\right)^2+((\tilde
u+c)H_{ij}+F_{ij})b^ib^j+G_ib^i+K.\end{align} Equation \eqref{tGk}
implies the existence of  $b^j(\tilde u)$ such that $\tilde
G_k=0$. Using the last equation, we can chose $d(\tilde u)$ such
that $\tilde K=0$. Equation \eqref{tHkl} implies the existence of
an orthogonal matrix $a_i^j$ such that $\tilde H_{kl}$ is a
diagonal matrix with the diagonal elements
$\lambda_1,...,\lambda_n$ such that $\lambda_1\leq \cdots\leq
\lambda_n$.

Since $\nabla R\neq 0$,  Corollary \ref{Corsym} shows that
$H_{ij}$ is not zero.

The transformation \eqref{trans1}  does not change the form of the
metric from Theorem \ref{ThMain} if and only if
$H_{kl}a^k_ia^l_j=H_{ij}$ and $b^i(\tilde u)$, $d(\tilde u)$
satisfy certain conditions. This and \eqref{tFkl} prove the last
claim of the theorem. \qed

{\bf Acknowledgements.} The authors would like to express their
thanks to the referee for his various comments that helped to
improve the style of the paper and to make the statement of the
main theorem more precise. The first author thanks Hamburg
University for support during preparation of this paper. The
second author is supported by the grant 201/09/P039 of the Grant
Agency of Czech Republic and by the grant MSM~0021622409 of the
Czech Ministry of Education.

\end{document}